\documentstyle[12pt]{article}
\newtheorem{thm}{Theorem}[section]
\newtheorem{Theorem}{Theorem}[section]

\newcommand{\sthm}{\begin{Theorem}}         
\newcommand{\ethm}{\end{Theorem}}           
\newtheorem{lem}[thm]{Lemma}
\newtheorem{Corollary}[Theorem]{Corollary}   
\newcommand{\scor}{\begin{Corollary}}       
\newcommand{\ecor}{\end{Corollary}}         
\newtheorem{defn}{{\sc Definition}}[section]
\newcommand{\pf}{ \par \vspace{1ex} \noindent {\sc Proof.} \hspace{2mm}}

\begin{document}
\title{The existence, nonexistence and uniqueness of global positive coexistence of a nonlinear elliptic biological interacting model}
\author{Joon Hyuk Kang, Yun Myung Oh\\
Department of Mathematics, Andrews University\\
Berrien Springs, MI 49104, U.S.A.\\
(kang@andrews.edu, ohy@andrews.edu)}
\date{ }
\maketitle
\begin{abstract}
The purpose of this paper is to give a sufficient condition for the existence, nonexistence and uniqueness of coexistence of positive solutions to a rather general type of elliptic competition system of the Dirichlet problem on the bounded domain $\Omega$ in $R^{n}$. The techniques used in this paper are upper-lower solutions, maximum principles and spectrum estimates. The arguments also rely on some detailed properties for the solution of logistic equations. This result yields an algebraically computable criterion for the positive coexistence of competing species of animals in many biological models.
\\[0.2in]
Keywords:competition model, coexistence state
\\[0.2in]
AMS 2000 Mathematics Subject Classification: 35A05, 35A07, 35B50, 35G30, 35J25 and 35K20
\\[0.2in]
Research supported by Andrews University Faculty Research Grant 2000 
\end{abstract}
\section{Introduction}
The coexistence of steady states of competition interacting models with diffusion has been an object of intensive study in recent years. See, for example, lists of references in ~\cite{ac92}, ~\cite{cc87}, ~\cite{cc89}, ~\cite{cl84}, ~\cite{gl94}, ~\cite{gp91}, ~\cite{ll91}, ~\cite{l80}, ~\cite{ll91}, ~\cite{l80}.   
The most general type of parabolic competition interacting system is 
$$\left\{\begin{array}{l}
u_{t} = \Delta u + ug(u,v),\\
v_{t} = \Delta v + vh(u,v),
\end{array}\right.$$
where $\Delta$ is the Laplacian and $u, v$ represent the densities of two competing species of animals. The terms $\Delta u$ and $\Delta v$ model dispersal by means of simple diffusion. We assume here that the $C^{1}$ functions $g$ and $h$ are relative growth rates satisfying the following so-called growth rate conditions:
\\[0.1in]
$(G_{1})$ $g_{u}(u,v) < 0, g_{v}(u,v) < 0, h_{u}(u,v) < 0, h_{v}(u,v) < 0$,\\
$(G_{2})$ There exist constants $c_{0} > 0, c_{1} > 0$ such that $g(u,0) \leq 0$ for $u \geq c_{0}$ and $h(0,v) \leq 0$ for $u \geq c_{1}$.
\vskip1pc
The hypothesis $(G_{1})$ characterizes how the two species $u$ and $v$ interact with each other in terms of their relative growth rates. It is well known that condition $(G_{2})$ exhibits the so-called logistic pattern while the constants $c_{0}$ and $c_{1}$ are referred to as the carrying capacity.
\vskip1pc
The earlier literature on this line focused in the Neumann boundary value problem:
\begin{equation}
\left\{\begin{array}{l}
u_{t} = \Delta u + ug(u,v),\\
v_{t} = \Delta v + vh(u,v),\\
\frac{\partial u(t,x)}{\partial n} = 0 = \frac{\partial v(t,x)}{\partial n}\;\;\mbox{for}\;\;(t,x) \in [0,T] \times \partial\Omega,
\end{array}\right.\label{eq:1}
\end{equation}  
and in its steady state, the elliptic system
\begin{equation}
\left\{\begin{array}{l}
\Delta u + ug(u,v) = 0,\\
\Delta v + vh(u,v) = 0,\\
\frac{\partial u(x)}{\partial n} = 0 = \frac{\partial v(x)}{\partial n}\;\;\mbox{for}\;\;x \in \partial\Omega,
\end{array}\right.\label{eq:2}
\end{equation}  
where $n$ denotes the unit out-normal along boundary $\partial\Omega$. The Neumann boundary conditions $\frac{\partial u(x)}{\partial n} = 0 = \frac{\partial v(x)}{\partial n}$ are interpreted as an assumption that both populations are staying inside, that there is no migratory flux across $\partial\Omega$. The goals of investigations along this line include finding out under what conditions on the nonlinearities $g$ and $h$ systems (\ref{eq:1}) and (\ref{eq:2}) have positive solutions $u > 0, v > 0$ and the possible uniqueness. Most of the work in this case were established by P. DeMottoni and F. Rothe in 1979 ~\cite{dr79} and P. Brown in 1980 ~\cite{b80}. Their work in a large sense completes the avenue of investigation in the study of Neumann boundary value problems. Researchers thus have since turned their attention to the biologically and physically more important case that is the Dirichlet boundary condition:
\begin{equation}
\left\{\begin{array}{rcl}
\Delta u + ug(u,v) & = & 0,\\
\Delta v + vh(u,v) & = & 0,\\
(u,v)|_{\partial\Omega} & = & (0,0).
\end{array}\right.\label{eq:3}
\end{equation}
Biologically, this setting allows migration of these two populations across the boundary but they may not stay on $\partial\Omega$, where, for example $\partial\Omega$ is a river. It was then found that the features known in the Neumann setting are not usually shared by those in the Dirichlet setting. The study in the latter setting, especially in the case of steady states like system (\ref{eq:3}), seems to be more difficult.
\vskip1pc
The goal of this paper is to answer the following questions about positive steady state to (\ref{eq:3}).  
\\[0.2in]
\underline{Problem 1} : What are the sufficient conditions for existence of steady state? 
\\[0.2in]
\underline{Problem 2} : Is it possible for either one of the species to be extinct?
\\[0.2in]
\underline{Problem 3} : When is the coexistence state unique? 
\\[0.2in]

\section{Preliminaries}
In this section we state some preliminary results which will be useful for 
our later arguments.
\begin{defn}\label{defn:2.1}(Upper and Lower solutions)\\
The vector functions $(\bar{u}^{1},...,\bar{u}^{N}), (\underline{u}^{1},...,\underline{u}^{N})$ form an upper/lower solution pair for the system
$$\left\{\begin{array}{r}
\Delta u^{i} + g^{i}(u^{1},...,u^{N}) = 0\;\;\mbox{in}\;\;\Omega\\
u^{i} = 0\;\;\mbox{on}\;\;\partial\Omega
\end{array}\right.$$
if for $i = 1,...,N$
$$\left\{\begin{array}{r}
\Delta \bar{u}^{i} + g^{i}(u^{1},...,u^{i-1},\bar{u}^{i},u^{i+1},...,u^{N}) \leq 0\\
\Delta \underline{u}^{i} + g^{i}(u^{1},...,u^{i-1},\underline{u}^{i},u^{i+1},...,u^{N}) \geq 0\\
\mbox{in}\;\;\Omega\;\;\mbox{for}\;\;\underline{u}^{j} \leq u^{j} \leq \bar{u}^{j}, j \neq i,
\end{array}\right.$$
and
$$\left.\begin{array}{l}
\underline{u}^{i} \leq \bar{u}^{i}\;\;\mbox{on}\;\;\Omega\\
\underline{u}^{i} \leq 0 \leq \bar{u}^{i}\;\;\mbox{on}\;\;\partial\Omega.
\end{array}\right.$$
\end{defn}
\begin{lem}\label{lem:2.1}(~\cite{ac92})\\
If $g^{i}$ in the Definition \ref{defn:2.1} are in $C^{1}$ and the system admits an upper/lower solution pair $(\underline{u}^{1},...,\underline{u}^{N}), (\bar{u}^{1},...,\bar{u}^{N})$, then there is a solution of the system in \ref{defn:2.1} with $\underline{u}^{i} \leq u^{i} \leq \bar{u}^{i}$ in $\bar{\Omega}$. If
$$\left.\begin{array}{l}
\Delta\bar{u}^{i} + g^{i}(\bar{u}^{1},...,\bar{u}^{N}) \neq 0,\\
\Delta\underline{u}^{i} + g^{i}(\underline{u}^{1},...,\underline{u}^{N}) \neq 0
\end{array}\right.$$
in $\Omega$ for $i = 1,...,N$, then $\underline{u}^{i} < u^{i} < \bar{u}^{i}$ in $\Omega$.
\end{lem}

\begin{lem}\label{lem:2.2}(The first eigenvalue)(\cite{d})
\begin{equation}
\left\{ \begin{array}{l}
-\Delta u + q(x)u = \lambda u\;\;\mbox{in}\;\;\Omega,\\
u|_{\partial \Omega} = 0,
\end{array} \right. \label{eq:4}
\end{equation}
where $q(x)$ is a smooth function from $\Omega$ to $R$ and $\Omega$ is a bounded domain in $R^{n}$.\\
$(A)$ The first eigenvalue $\lambda_{1}(q)$ of (\ref{eq:4}), denoted by simply $\lambda_{1}$ when $q \equiv 0$, is simple with a positive 
eigenfunction.\\
$(B)$ If $q_{1}(x) < q_{2}(x)$ for all $x \in \Omega$, then $\lambda_{1}(q_{1}) 
< \lambda_{1}(q_{2})$.\\
$(C)$(Variational Characterization of the first eigenvalue)\\
$$\lambda_{1}(q) = \min_{\phi \in W_{0}^{1}(\Omega),\phi \neq 
0}\frac{\int_{\Omega}(|\nabla \phi|^{2}+q\phi^{2})dx}{\int_{\Omega}\phi^{2}dx}$$.
\end{lem}
\vskip1pc
We also need some information on the solutions of the following logistic equations.
\begin{lem}\label{lem:2.4}(\cite{ll91})
$$\left\{ \begin{array}{l}
\Delta u + uf(u) = 0\;\; \mbox{in}\;\; \Omega,\\
u|_{\partial\Omega} = 0, u > 0,
\end{array} \right.$$ 
where $f$ is a decreasing $C^{1}$ function such that there exists $c_{0} > 0$ 
such that $f(u) \leq 0$ for $u \geq c_{0}$ and $\Omega$ is a bounded domain in $R^{n}$.
\\[0.1in]
$(1)$ If $f(0) > \lambda_{1}$, then the 
above equation has a unique positive solution, where 
$\lambda_{1}$ is the first eigenvalue of $-\Delta$ with homogeneous boundary 
condition. We denote this unique positive solution as $\theta_{f}$.\\
$(2)$ If $f(0) \leq \lambda_{1}$, then the above equation does not have any positive solution.
\end{lem}

\section{Existence, Nonexistence and Uniqueness}
We consider the system (\ref{eq:3}) with conditions $(G_{1})$ and $(G_{2})$.
\begin{thm}\label{thm:3.1}
$(A)$ If $g(0,c_{1}) > \lambda_{1}$ and $h(c_{0},0) > \lambda_{1}$, then (\ref{eq:3}) has a solution $(u,v)$ with
$$\left.\begin{array}{l}
\theta_{g(\cdot,c_{1})} < u < \theta_{g(\cdot,0)}\\
\theta_{h(c_{0},\cdot)} < v < \theta_{h(0,\cdot)}.
\end{array}\right.$$
Conversely, any solution $(u,v)$ of (\ref{eq:3}) with $u > 0, v > 0$ in $\Omega$ must satisfy these inequalities.\\
$(B)$ If $g(0,0) \leq \lambda_{1}$ or $h(0,0) \leq \lambda_{1}$, then (\ref{eq:3}) does not have any positive solution.
\end{thm}
\pf\\
$(A)$ Let $\bar{u} = \theta_{g(\cdot,0)}, \bar{v} = \theta_{h(0,\cdot)}$. Then by the monotonicity of $g$,
$$\left.\begin{array}{ll}
 & \Delta\bar{u} + \bar{u}g(\bar{u},\bar{v})\\
= & \Delta\bar{u} + \bar{u}(g(\bar{u},0) - g(\bar{u},0) + g(\bar{u},\bar{v}))\\
= & \bar{u}(g(\bar{u},\bar{v}) - g(\bar{u},0)) < 0.
\end{array}\right.$$
Similarly, $$\Delta\bar{v} + \bar{v}h(\bar{u},\bar{v}) < 0.$$
So, $(\bar{u},\bar{v})$ is an upper solution to (\ref{eq:3}).\\
Let $\underline{u} = \theta_{g(\cdot,c_{1})}$ and $\underline{v} = \theta_{h(c_{0},\cdot)}$. Then by the Maximum Principles, we obtain
$$\left.\begin{array}{l}
\underline{u} \leq \theta_{g(\cdot,0)} \leq c_{0},\\
\underline{v} \leq \theta_{h(0,\cdot)} \leq c_{1}.
\end{array}\right.$$
By the monotonicity of $g$,
$$\left.\begin{array}{ll}
 & \Delta\underline{u} + \underline{u}g(\underline{u},\underline{v})\\
= & \Delta\underline{u} + \underline{u}(g(\underline{u},c_{1}) - g(\underline{u},c_{1}) + g(\underline{u},\underline{v}))\\
= & \underline{u}(g(\underline{u},\underline{v}) - g(\underline{u},c_{1})) \geq 0.
\end{array}\right.$$
Similarly, $$\Delta\underline{v} + \underline{v}h(\underline{u},\underline{v}) \geq 0.$$
Therefore, $(\underline{u},\underline{v})$ is a lower solution to (\ref{eq:3}). Furthermore, $\underline{u} < \bar{u}, \underline{v} < \bar{v}$ in $\Omega$ and $\underline{u} = \bar{u} = \underline{v} = \bar{v} = 0$ on $\partial\Omega$.\\
So, (\ref{eq:3}) has a solution $(u,v)$ with
$$\left.\begin{array}{l}
\theta_{g(\cdot,c_{1})} < u < \theta_{g(\cdot,0)},\\
\theta_{h(c_{0},\cdot)} < v < \theta_{h(0,\cdot)}.
\end{array}\right.$$
Suppose $(u,v)$ is a coexistence state for (\ref{eq:3}). Then since
$$\left.\begin{array}{ll}
 & \Delta u + ug(u,0)\\
\geq & \Delta u + ug(u,v) = 0,
\end{array}\right.$$
$u$ is a lower solution of 
\begin{equation}
\left.\begin{array}{rrl}
\Delta Z + Zg(Z,0) & = & 0\;\;\mbox{in}\;\;\Omega,\\
Z & = & 0\;\;\mbox{on}\;\;\partial\Omega.
\end{array}\right.\label{eq:5}
\end{equation}
But, since any constant larger than $c_{0}$ is an upper solution of (\ref{eq:5}), we have
\begin{equation}
\left.\begin{array}{l}
u < \theta_{g(\cdot,0)}.
\end{array}\right.\label{eq:6}
\end{equation}
Similarly, we have
\begin{equation}
\left.\begin{array}{l}
v < \theta_{h(0,\cdot)}.
\end{array}\right.\label{eq:6'}
\end{equation}
Since $v < \theta_{h(0,\cdot)} \leq c_{1}$, by the monotonicity of $g$
$$\left.\begin{array}{ll}
 & \Delta u + ug(u,c_{1})\\
\leq & \Delta u + ug(u,v) = 0.
\end{array}\right.$$
Therefore, $u$ is an upper solution of 
\begin{equation}
\left.\begin{array}{rrl}
\Delta Z + Zg(Z,c_{1}) & = & 0\;\;\mbox{in}\;\;\Omega,\\
Z & = & 0\;\;\mbox{on}\;\;\partial\Omega.
\end{array}\right.\label{eq:7}
\end{equation}
If $\epsilon > 0$ is so small that $g(\epsilon\phi_{1},c_{1}) > \lambda_{1}$ on $\bar{\Omega}$, where $\phi_{1}$ is the first eigenvector of $-\Delta$ with homogeneous boundary condition, then since
$$\left.\begin{array}{ll}
 & \Delta\epsilon\phi_{1} + \epsilon\phi_{1}g(\epsilon\phi_{1},c_{1})\\
= & \epsilon(\Delta\phi_{1} + \phi_{1}g(\epsilon\phi_{1},c_{1}))\\
> & \epsilon(\Delta\phi_{1} + \lambda_{1}\phi_{1}) = 0,
\end{array}\right.$$
$\epsilon\phi_{1}$ is a lower solution of (\ref{eq:7}). So, we have
\begin{equation}
\left.\begin{array}{l}
\theta_{g(\cdot,c_{1})} < u.
\end{array}\right.\label{eq:8}
\end{equation}
Similarly, we have 
\begin{equation}
\left.\begin{array}{l}
\theta_{h(c_{0},\cdot)} < v.
\end{array}\right.\label{eq:8'}
\end{equation}
By (\ref{eq:6}), (\ref{eq:6'}), (\ref{eq:8}) and (\ref{eq:8'}),
$$\left.\begin{array}{l}
\theta_{g(\cdot,c_{1})} < u < \theta_{g(\cdot,0)},\\
\theta_{h(c_{0},\cdot)} < v < \theta_{h(0,\cdot)}.
\end{array}\right.$$
\\[0.1in]
$(B)$ Assume $g(0,0) \leq \lambda_{1}$. The other cases are proved similarly. Suppose $(\bar{u},\bar{v})$ is a positive solution to (\ref{eq:3}). Then since
$$\left.\begin{array}{ll}
 & \Delta\bar{u} + \bar{u}g(\bar{u},0)\\
= & \Delta\bar{u} + \bar{u}(g(\bar{u},\bar{v}) - g(\bar{u},\bar{v}) + g(\bar{u},0))\\
= & \bar{u}(g(\bar{u},0) - g(\bar{u},\bar{v})) \geq 0,
\end{array}\right.$$
$\bar{u}$ is a lower solution to 
\begin{equation}
\left.\begin{array}{rrl}
\Delta u + ug(u,0) & = & 0\;\;\mbox{in}\;\;\Omega,\\
u & = & 0\;\;\mbox{on}\;\;\partial\Omega.
\end{array}\right.\label{eq:9}
\end{equation}
Any constant larger than $c_{0}$ is an upper solution to (\ref{eq:9}). Hence, (\ref{eq:9}) has a positive solution $u_{0}$ with $\bar{u} < u_{0}$. This contradicts to the Lemma \ref{lem:2.4} which says there is no positive solution of (\ref{eq:9}) if $g(0,0) \leq \lambda_{1}$.

\begin{thm}\label{thm:3.2} If $g(0,c_{1}) > \lambda_{1}, h(c_{0},0) > \lambda_{1}$ and
$$\left.\begin{array}{lll}
4\inf(-\frac{\partial g(u,v)}{\partial u})\inf(-\frac{\partial h(u,v)}{\partial v}) & \geq & \frac{\theta_{g(\cdot,0)}}{\theta_{h(c_{0},\cdot)}}(\sup\frac{\partial g(u,v)}{\partial v})^{2} + \frac{\theta_{h(0,\cdot)}}{\theta_{g(\cdot,c_{1})}}(\sup\frac{\partial h(u,v)}{\partial u})^{2}\\
 & & + 2(\sup\frac{\partial g(u,v)}{\partial v})(\sup\frac{\partial h(u,v)}{\partial u}),
\end{array}\right.$$ 
then (\ref{eq:3}) has a unique positive solution. 
\end{thm}
\pf Suppose $(u_{1},v_{1})$ and $(u_{2},v_{2})$ are positive solutions to (\ref{eq:3}). Let $p = u_{1} - u_{2}$ and $q = v_{1} - v_{2}$. Then
$$\left.\begin{array}{ll}
 & \Delta p + pg(u_{1},v_{1})\\
= & \Delta u_{1} - \Delta u_{2} + (u_{1} - u_{2})g(u_{1},v_{1})\\
= & -\Delta u_{2} - u_{2}g(u_{1},v_{1})\\
= & -\Delta u_{2} - u_{2}(g(u_{2},v_{2}) - g(u_{2},v_{2}) + g(u_{1},v_{1}))\\
= & -u_{2}(g(u_{1},v_{1}) - g(u_{2},v_{2}))\\
= & -u_{2}(g(u_{1},v_{1}) - g(u_{2},v_{1}) + g(u_{2},v_{1}) - g(u_{2},v_{2})).
\end{array}\right.$$
But, by the Mean Value Theorem, there is $\tilde{x}$ depending on $u_{1}, u_{2}$ such that
$$g(u_{1},v_{1}) - g(u_{2},v_{1}) = \frac{\partial g(\tilde{x},v_{1})}{\partial u}p.$$
Hence,
$$\Delta p + pg(u_{1},v_{1}) = -u_{2}[\frac{\partial g(\tilde{x},v_{1})}{\partial u}p + g(u_{2},v_{1}) - g(u_{2},v_{2})].$$
i.e.,\\
\begin{equation}
\left.\begin{array}{l}
\Delta p + g(u_{1},v_{1})p + u_{2}p\frac{\partial g(\tilde{x},v_{1})}{\partial u} - u_{2}(g(u_{2},v_{2}) - g(u_{2},v_{1})) = 0.
\end{array}\right.\label{eq:10}
\end{equation}
The same argument shows that 
\begin{equation}
\left.\begin{array}{l}
\Delta q + h(u_{2},v_{2})q + v_{1}q\frac{\partial h(u_{2},\bar{x})}{\partial v} - v_{1}(h(u_{2},v_{1}) - h(u_{1},v_{1})) = 0,
\end{array}\right.\label{eq:11}
\end{equation}
where $\bar{x}$ depends on $v_{1}, v_{2}$ by the Mean Value Theorem.\\
Since $\lambda_{1}(-g(u_{1},v_{1})) = 0$, by the Variational Characterization of the first eigenvalue,
\begin{equation}
\left.\begin{array}{l}
\int_{\Omega}Z(-\Delta Z - g(u_{1},v_{1})Z)dx \geq 0
\end{array}\right.\label{eq:12}
\end{equation}
for any $Z \in C^{2}(\bar{\Omega})$ and $Z|_{\partial\Omega} = 0$. The same argument shows that 
\begin{equation}
\left.\begin{array}{l}
\int_{\Omega}W(-\Delta W - h(u_{2},v_{2})W)dx \geq 0
\end{array}\right.\label{eq:13}
\end{equation}
for any $W \in C^{2}(\bar{\Omega})$ and $W|_{\partial\Omega} = 0$.\\
From (\ref{eq:10}) and (\ref{eq:11}), we get
$$\left.\begin{array}{l}
-p\Delta p - g(u_{1},v_{1})p^{2} - \frac{\partial g(\tilde{x},v_{1})}{\partial u}u_{2}p^{2} + u_{2}p(g(u_{2},v_{2}) - g(u_{2},v_{1})) = 0,\\
-q\Delta q - h(u_{2},v_{2})q^{2} - \frac{\partial h(u_{2},\bar{x})}{\partial v}v_{1}q^{2} + v_{1}q(h(u_{2},v_{1}) - h(u_{1},v_{1})) = 0.
\end{array}\right.$$
Hence from (\ref{eq:12}) and (\ref{eq:13}),
$$\left.\begin{array}{l}
\int_{\Omega}(-\frac{\partial g(\tilde{x},v_{1})}{\partial u}u_{2}p^{2} + u_{2}p(g(u_{2},v_{2}) - g(u_{2},v_{1})) + v_{1}q(h(u_{2},v_{1}) - h(u_{1},v_{1}))\\
- \frac{\partial h(u_{2},\bar{x})}{\partial v}v_{1}q^{2})dx \leq 0.
\end{array}\right.$$
By the Mean Value Theorem, for each $x \in \Omega$, there exist $\tilde{y}, \bar{y}$ such that
$$\left.\begin{array}{l}
g(u_{2},v_{2}) - g(u_{2},v_{1}) = \frac{\partial g(u_{2},\tilde{y})}{\partial v}(-q),\\
h(u_{2},v_{1}) - h(u_{1},v_{1}) = \frac{\partial h(\bar{y},v_{1})}{\partial u}(-p),
\end{array}\right.$$
which implies that
$$\int_{\Omega}-\frac{\partial g(\tilde{x},v_{1})}{\partial u}u_{2}p^{2} - (u_{2}\frac{\partial g(u_{2},\tilde{y})}{\partial v} + v_{1}\frac{\partial h(\bar{y},v_{1})}{\partial u})pq - \frac{\partial h(u_{2},\bar{x})}{\partial v}v_{1}q^{2}dx \leq 0.$$
Therefore, we find
\begin{eqnarray*}
p \equiv q \equiv 0 & \mbox{if} & -\frac{\partial g(\tilde{x},v_{1})}{\partial u}u_{2}\zeta^{2} - 
(u_{2}\frac{\partial g(u_{2},\tilde{y})}{\partial v} + v_{1}\frac{\partial h(\bar{y},v_{1})}{\partial u})\zeta\eta\\
& - &
\frac{\partial h(u_{2},\bar{x})}{\partial v}v_{1}\eta^{2}
\;\; \mbox{is 
positive definite}\\
& & \mbox{for each $x \in \Omega.$}
\end{eqnarray*}                    
This is the case if
$$\left.\begin{array}{l}
u_{2}^{2}(\frac{\partial g(u_{2},\tilde{y})}{\partial v})^{2} + 
v_{1}^{2}(\frac{\partial h(\bar{y},v_{1})}{\partial u})^{2} + 
2u_{2}v_{1}\frac{\partial g(u_{2},\tilde{y})}{\partial v}\frac{\partial h(\bar{y},v_{1})}{\partial u}\\
- 4\frac{\partial g(\tilde{x},v_{1})}{\partial u}\frac{\partial h(u_{2},\bar{x})}{\partial v}u_{2}v_{1} \leq 0 \;\;\mbox{for each $x 
\in \Omega.$}
\end{array}\right.$$
$$\left.\begin{array}{lll}
\mbox{i.e.,}\;\;4\frac{\partial g(\tilde{x},v_{1})}{\partial u}\frac{\partial h(u_{2},\bar{x})}{\partial v} & 
\geq & 
\frac{u_{2}}{v_{1}}(\frac{\partial g(u_{2},\tilde{y})}{\partial v})^{2} + 
\frac{v_{1}}{u_{2}}(\frac{\partial h(\bar{y},v_{1})}{\partial u})^{2}\\
 & & + 
2\frac{\partial g(u_{2},\tilde{y})}{\partial v}\frac{\partial h(\bar{y},v_{1})}{\partial u}\;\; \mbox{for each $x \in \Omega$}.
\end{array}\right.$$ 
But, from the inequality in $(A)$ and the hypothesis in the theorem,
\begin{eqnarray*}
\left. \begin{array}{lll}
\frac{u_{2}}{v_{1}}(\frac{\partial g(u_{2},\tilde{y})}{\partial v})^{2} &+& 2\frac{\partial g(u_{2},\tilde{y})}{\partial v}\frac{\partial h(\bar{y},v_{1})}{\partial u} 
\;\;+\;\; 
\frac{v_{1}}{u_{2}}(\frac{\partial h(\bar{y},v_{1})}{\partial u})^{2}\\
& \leq & 
\frac{\theta_{g(\cdot,0)}}{\theta_{h_(c_{0},\cdot)}}(\sup\frac{\partial g(u,v)}{\partial v})^{2} + 
\frac{\theta_{h(0,\cdot)}}{\theta_{g(\cdot,c_{1})}}(\sup\frac{\partial h(u,v)}{\partial u})^{2}\\
& & + 2\sup(\frac{\partial g(u,v)}{\partial v})\sup(\frac{\partial h(u,v)}{\partial u})\\
& \leq & 4\inf(-\frac{\partial g(u,v)}{\partial u})\inf(-\frac{\partial h(u,v)}{\partial v})\\
& \leq & 4\frac{\partial g(\tilde{x},v_{1})}{\partial u}\frac{\partial g(u_{2},\bar{x})}{\partial v}.
\end{array} \right. 
\end{eqnarray*}
We can also extend the results to the case when there are multiple species competing in the same environment.
\\[0.1in]
Consider the interacting model
\begin{equation}
\left.\begin{array}{rrl}
\Delta u_{i} + u_{i}g_{i}(u_{i},u_{2},...,u_{i},u_{i+1},...,u_{N}) & = & 0\;\;\mbox{in}\;\;\Omega,\\
u_{i} & = & 0\;\;\mbox{on}\;\;\partial\Omega
\end{array}\right.\label{eq:16}
\end{equation}
for $i = 1,...,N$.
\\[0.1in]
Again, we assume here that the $C^{1}$ functions $g_{i}$ for $i = 1,...,N$ are relative growth rates satisfying the following growth rate conditions:
\\[0.1in]
$(M1)$ $\frac{\partial g_{i}}{\partial u_{j}} < 0$ for $i, j = 1,2,...,N$,\\
$(M2)$ There exist constants $c_{1} > 0, c_{2} > 0,...,c_{N} > 0$ such that\\
$g_{i}(0,...,0,u_{i},0,...,0) \leq 0$ for $u_{i} \geq c_{i}$.
\\[0.1in]
Again, $(M1)$ characterizes how the $N$ species $u_{1}, u_{2},...,u_{N}$ interact with each other in terms of their relative growth rates and $(M2)$ is the logistic pattern with carrying capacity constants $c_{1}, c_{2},...,c_{N}$.
\\[0.in]
The followings are the main results. The proofs are similar to those with 2 competing species, and so we just sketch it without the details.
\begin{thm}\label{thm:4.1} 
$(A)$ If $g_{i}(c_{1},c_{2},...,c_{i-1},0,c_{i+1},...,c_{N}) > \lambda_{1}$ for $i = 1,...,N$, then (\ref{eq:16}) has a solution $(u_{1},...,u_{N})$ with
$$\theta_{g_{i}(c_{1},...,c_{i-1},\cdot,c_{i+1},...,c_{N})} < u_{i} < \theta_{g_{i}(0,...,0,\cdot,0,...,0)}$$
for $i = 1,..,N$.\\
Conversely, any solution $(u_{1},...,u_{N})$ of (\ref{eq:16}) with $u_{i} > 0$ in $\Omega$ must satisfy these inequalities.\\
$(B)$ If $g_{i}(0,...,0) \leq \lambda_{1}$ for some $i = 1,...,N$, then (\ref{eq:16}) does not have any positive solution.
\end{thm}
\pf\\
$(A)$ Let $\bar{u_{i}} = \theta_{g_{i}(0,...,0,\cdot,0,...,0)}$ and $\underline{u_{i}} = \theta_{g_{i}(c_{1},...,c_{i-1},\cdot,c_{i+1},...,c_{N})}$  for $i = 1,...,N$. Then by the Maximum Principles and the monotonicity of $g_{i}$, $(\bar{u_{1}},...,\bar{u_{i}},...,\bar{u_{N}})$ and $(\underline{u_{1}},...,\underline{u_{i}},...,\underline{u_{N}})$ are upper and lower solutions to (\ref{eq:16}), respectively.\\
Furthermore, for $i = 1,...,N$, $\underline{u_{i}} < \bar{u_{i}}$ in $\Omega$ and $\underline{u_{i}} = \bar{u_{i}} = 0$ on $\partial\Omega$.\\
So, (\ref{eq:16}) has a solution $(u_{1},...,u_{N})$ with the desired inequalities
$$\theta_{g_{i}(c_{1},...,c_{i-1},\cdot,c_{i+1},...,c_{N})} < u_{i} < \theta_{g_{i}(0,...,0,\cdot,0,...,0)}$$
for $i = 1,..,N$.\\
Suppose $(u_{1},...,u_{N})$ is a coexistence state for (\ref{eq:16}). Then by the direct computation using the monotonicity of $g_{i}$, we know that $u_{i}$ is a lower solution of
\begin{equation}
\left.\begin{array}{rrl}
\Delta Z + Zg_{i}(0,...,0,Z,0,...,0) & = & 0\;\;\mbox{in}\;\;\Omega,\\
Z & = & 0\;\;\mbox{on}\;\;\partial\Omega
\end{array}\right.\label{eq:17}
\end{equation}
for $i = 1,...,N$\\
But, since any constant larger than $c_{i}$ is an upper solution of (\ref{eq:17}), we have
\begin{equation}
\left.\begin{array}{l}
u_{i} < \theta_{g_{i}(0,...,0,\cdot,0,...,0)}
\end{array}\right.\label{eq:18}
\end{equation}
for $i = 1,...,N$.\\
Since $u_{i} < \theta_{g_{i}(0,...,0,\cdot,0,...,0)} \leq c_{i}$, by the monotonicity of $g_{i}$, we can derive that $u_{i}$ is an upper solution of
\begin{equation}
\left.\begin{array}{rrl}
\Delta Z + Zg_{i}(c_{1},...,c_{i-1},Z,c_{i+1},...,c_{N}) & = & 0\;\;\mbox{in}\;\;\Omega,\\
Z & = & 0\;\;\mbox{on}\;\;\partial\Omega
\end{array}\right.\label{eq:19}
\end{equation}
for $i = 1,...,N$.\\
If $\epsilon > 0$ is so small that $g_{i}(c_{1},...,c_{i-1},\epsilon\phi_{1},c_{i+1},...,c_{N}) > \lambda_{1}$ on $\bar{\Omega}$, where $\phi_{1}$ is the first eigenvector of $-\Delta$ with homogeneous boundary condition, then by the dirct computation again, we know that $\epsilon\phi_{1}$ is a lower solution of (\ref{eq:19}).\\
So, we have
\begin{equation}
\left.\begin{array}{l}
\theta_{g_{i}(c_{1},...,c_{i-1},\cdot,c_{i+1},...,c_{N})} < u_{i}
\end{array}\right.\label{eq:20}
\end{equation}
for $i = 1,...,N$.\\
By (\ref{eq:18}) and (\ref{eq:20}),
$$\theta_{g_{i}(c_{1},...,c_{i-1},\cdot,c_{i+1},...,c_{N})} < u_{i} < \theta_{g_{i}(0,...,0,\cdot,0,...,0)}$$
for $i = 1,...,N$.
\\[0.1in]
$(B)$ Without loss of generality, assume $g_{1}(0,...,0) \leq \lambda_{1}$.\\
Suppose $(\bar{u_{1}},...,\bar{u_{N}})$ is a positive solution to (\ref{eq:16}). Then by the monotonicity of $g_{i}$, $\bar{u_{1}}$ is a lower solution to
\begin{equation}
\left.\begin{array}{rrl}
\Delta Z + Zg_{1}(Z,0,...,0) & = & 0\;\;\mbox{in}\;\;\Omega,\\
Z & = & 0\;\;\mbox{on}\;\;\partial\Omega.
\end{array}\right.\label{eq:21}
\end{equation}
Hence, by the fact that any constant larger than $c_{1}$ is an upper solution to (\ref{eq:21}), (\ref{eq:21}) has a positive solution $u_{1}$ with $\bar{u_{1}} < u_{1}$ that contradicts to the Lemma \ref{lem:2.4}.

\begin{thm}\label{thm:4.2} If $g_{i}(c_{1},...,c_{i-1},0,c_{i+1},...,c_{N}) > \lambda_{1}$ and 
$$2\inf(-\frac{\partial g_{i}}{\partial x_{i}}) > \sum_{j=1,j\neq i}^{N}(\sup(-\frac{\partial g_{i}}{\partial x_{j}}) + K\sup(-\frac{\partial g_{j}}{\partial x_{i}}))$$
for $i = 1,...,N$, where $K = \sup_{i, j\neq i}\frac{\theta_{g_{j}(0,...,0,\cdot,0,...,0)}}{\theta_{g_{i}(c_{1},...,c_{i-1},\cdot,c_{i+1},...,c_{N})}},$ then (\ref{eq:16}) has a unique coexistence state.
\end{thm}
\pf Suppose $(u_{1},...,u_{N})$ and $(v_{1},...,v_{N})$ are coexistence states of (\ref{eq:16}) and let $w_{i} = u_{i} - v_{i}$ for $i = 1,...,N$. Then by the direct computation and the Variational Characterization of the first eigenvalu, we obtain
$$\int_{\Omega}\sum_{i=1}^{N}[v_{i}w_{i}(g_{i}(v_{1},...,v_{i},...,v_{N}) - g_{i}(u_{1},...,u_{i},...,u_{N}))]dx \leq 0.$$
By the Mean Value Theorem, there exist $t^{i}$ and $z^{ij}$ such that
\begin{equation}
\left.\begin{array}{l}
\int_{\Omega}\sum_{i=1}^{N}[\frac{\partial g_{i}(v_{1},...,v_{i-1},t^{i},v_{i+1},...,v_{N})}{\partial x_{i}}(-v_{i})w_{i}^{2}\\
+ \sum_{j=1,j\neq i}^{N}v_{i}w_{i}\frac{\partial g_{i}(u_{1},...,u_{j-1},z^{ij},v_{j+1},...,v_{N})}{\partial x_{j}}(- w_{j})]dx \leq 0.
\end{array}\right.\label{eq:24}
\end{equation}
If the integrand in the left side of (\ref{eq:24}) is positive definite, then (\ref{eq:24}) implies that $w_{i} \equiv 0$ in $\Omega$ for $i = 1,...,N$, which means the uniqueness of the coexistence state for (\ref{eq:16}). But for any $\epsilon > 0$,
$$\left.\begin{array}{ll}
 & \frac{\partial g_{i}(u_{1},...,u_{j-1},z^{ij},v_{j+1},...,v_{N})}{\partial x_{j}}(-v_{i})w_{i}(w_{j})\\
\leq & \frac{\partial g_{i}(u_{1},...,u_{j-1},z^{ij},v_{j+1},...,v_{N})}{\partial x_{j}}(-v_{i})[\frac{w_{i}^{2}}{2\epsilon} + \frac{\epsilon w_{j}^{2}}{2}].
\end{array}\right.$$
So, we can see that the integrand is positive definite if for $i = 1,...,N$ and $x \in \Omega$,
$$\left.\begin{array}{ll}
 & \frac{\partial g_{i}(v_{1},...,v_{i-1},t^{i},v_{i+1},...,v_{N})}{\partial x_{i}}(-v_{i})\\
> & \sum_{j=1,j\neq i}^{N}(\frac{\frac{\partial}{\partial x_{j}}g_{i}(u_{1},...,u_{j-1},z^{ij},v_{j+1},...,v_{N})(-v_{i})}{2\epsilon} + \frac{\epsilon\frac{\partial}{\partial x_{i}}g_{j}(u_{1},...,u_{i-1},z^{ji},v_{i+1},...,v_{N})(-v_{j})}{2})
\end{array}\right.$$
or equivalently,
\begin{equation}
\left.\begin{array}{ll}
& -\frac{\partial g_{i}(v_{1},...,v_{i-1},t^{i},v_{i+1},...,v_{N})}{\partial x_{i}}\\
> & \sum_{j=1,j\neq i}^{N}(\frac{-\frac{\partial}{\partial x_{j}}g_{i}(u_{1},...,u_{j-1},z^{ij},v_{j+1},...,v_{N})}{2\epsilon}\\
 & - \frac{\epsilon\frac{\partial}{\partial x_{i}}g_{j}(u_{1},...,u_{i-1},z^{ji},v_{i+1},...,v_{N})\frac{v_{j}}{v_{i}}}{2}).
\end{array}\right.\label{eq:25}
\end{equation}
Since $\theta_{g_{i}(c_{1},...,c_{i-1},\cdot,c_{i+1},...,c_{N})} < v_{i} < \theta_{g_{i}(0,...,0,\cdot,0,...,0)}$ in $\Omega$ for $i = 1,...,N$, (\ref{eq:25}) will hold if for $i = 1,...,N$,
$$\left.\begin{array}{ll}
& -\frac{\partial g_{i}(v_{1},...,v_{i-1},t^{i},v_{i+1},...,v_{N})}{\partial x_{i}}\\
> & \sum_{j=1,j\neq i}^{N}(\frac{\sup(-\frac{\partial g_{i}}{\partial x_{j}}}{2\epsilon} + \frac{\epsilon\sup(-\frac{\partial g_{j}}{\partial x_{i}})}{2}\frac{\theta_{g_{j}(0,...,0,\cdot,0,...,0)}}{\theta_{g_{i}(c_{1},...,c_{i-1},\cdot,c_{i+1},...,c_{N})}}).
\end{array}\right.$$
Let $K = \sup_{i, j\neq i}\frac{\theta_{g_{j}(0,...,0,\cdot,0,...,0)}}{\theta_{g_{i}(c_{1},...,c_{i-1},\cdot,c_{i+1},...,c_{N})}}.$ Then (\ref{eq:25}) holds if
$$\inf(-\frac{\partial g_{i}}{\partial x_{i}}) > \sum_{j=1,j\neq i}^{N}(\frac{\sup(-\frac{\partial g_{i}}{\partial x_{j}}}{2\epsilon} + \frac{K\epsilon\sup(-\frac{\partial g_{j}}{\partial x_{i}})}{2}.$$
Choosing $\epsilon = 1$, we have
$$2\inf(-\frac{\partial g_{i}}{\partial x_{i}}) > \sum_{j=1,j\neq i}^{N}(\sup(-\frac{\partial g_{i}}{\partial x_{j}}) + K\sup(-\frac{\partial g_{j}}{\partial x_{i}})).$$

\end{document}